\documentclass[preprint,12pt]{elsarticle}

\usepackage{amssymb}
\usepackage{graphicx}
\usepackage{mathptmx}      % use Times fonts if available on your TeX system

\usepackage{amssymb,amsfonts,amsmath,amscd,latexsym}
\usepackage[english]{babel}
\usepackage{geometry}
\usepackage{latexsym}
\usepackage{amssymb}
\usepackage[all]{xy}
\usepackage{amssymb}
\usepackage{latexsym}
\usepackage{graphicx}
\usepackage{appendix}
\usepackage[active]{srcltx}
\usepackage{srcltx}
\usepackage{xypic}
\usepackage{amsthm}

\newdefinition{defn}{Definition}[section]
\newdefinition{nim}[defn]{}
\newdefinition{rem}[defn]{Remark}
\newdefinition{ex}[defn]{Example}
\newtheorem{thm}{Theorem}[section]

\newtheorem{lem}[thm]{Lemma}

\begin{document}

\begin{frontmatter}



\title{Hochschild cohomology of fully group-graded algebras as Mackey functor}

\author{Tiberiu Cocone\c t}
\ead{tiberiu.coconet@math.ubbcluj.ro}
\address{Faculty of Economics and Business Administration, Babe\c s-Bolyai University, Str. Teodor Mihali, nr.58-60 , 400591 Cluj-Napoca, Romania}

\author{Constantin-Cosmin Todea\corref{cor1}}
\ead{Constantin.Todea@math.utcluj.ro}
\address{Department of Mathematics, Technical University of Cluj-Napoca, Str. G. Baritiu 25,
 Cluj-Napoca 400027, Romania}
\cortext[cor1]{Corresponding author} 

\begin{abstract} We prove that Hochschild cohomology of a certain class of fully group-graded algebras is a Mackey functor. We use the machinery of transfer maps between the Hochschild cohomology of symmetric algebras.
\end{abstract}


\end{frontmatter}

% \linenumbers

\section{Introduction }

Let $k$ be a field, $G$ be a finite group and let $R_G$ be a fully $G$-graded $k$-algebra. By definition, $R_G$ has a decomposition $R_G=\oplus_{g\in G} R_g,$ where each $R_g$ is a a $k$-vector space, for any $g\in G$, such that for all $g,h\in G$ we have $R_gR_h=R_{gh}.$  In this paper we work only with symmetric fully $G$-graded algebras $R_G$ with the property that for any  subgroup $H$ of $G$ we have that $R_H$ is a parabolic subalgebra of $R_G$, for further details see \cite[Definition 2.3, Definition 5.1]{Coconet-Todea-BrHig}.
We begin by giving an example of a group-graded algebra which lies in the above mentioned class of symmetric algebras and which includes the group algebra case.
\begin{ex}\label{Coconet-Todea-ex1} Let $R_G$ be a fully $G$-graded algebra such that $R_1$ is a full matrix ring. By \cite[Theorem 4]{Coconet-Todea-Har} we know that $R_G$ is a symmetric $k$-algebra. Moreover by \cite[Lemma 1.1]{Coconet-Todea-Schm} we obtain that $R_G$ is a crossed product. Furthermore, these statements assure us that for any  subgroup $H$ of $G$, $R_H$ is a parabolic subalgebra of $R_G$.
\end{ex}
We use the language of transfer maps between Hochschild cohomology algebras of symmetric algebras, defined originally in \cite{Coconet-Todea-LiTr}, to give a structure of Mackey functor for Hochschild cohomology of such fully $G$-graded algebras. Recall from \cite{Coconet-Todea-LiTr} that if $M$ is an $A-B$-bimodule, projective as a left $A$-module and as a right $B$-module (where $A,B$ are two symmetric $k$-algebras) there is a graded $k$-linear map, called transfer map
$$t_M:\mathrm{HH}^*(B)\longrightarrow\mathrm{HH}^*(A).$$
These transfer maps are also analyzed by explicit definitions in \cite{Coconet-Todea-LiZh}.

Let $H$ be a subgroup of $G$ and let $g\in G$. We denote by $M:=R_G$ the $R_H-R_G$-bimodule structure on $R_G$ given by multiplication in $R_G$; similarly $N:=R_G$ as an $R_G-R_H$-bimodule. Also we consider $P:=R_{[gH]}$ to be the $R_{[^gH]}-R_H$-bimodule given by multiplication in $R_G$.  By Lemma \ref{Coconet-Todea-lema}, a) (which is proved in the next section)  we can define
$$r_H^G=t_{M},~~t_H^G=t_{N},~~c_{g,H}=t_{P}, $$
hence we have the following graded $k$-linear maps
$$r_H^G:\mathrm{HH}^*(R_G)\rightarrow\mathrm{HH}^*(R_H);$$
$$t_H^G:\mathrm{HH}^*(R_H)\rightarrow \mathrm{HH}^*(R_G);$$
$$c_{g,H}:\mathrm{HH}^*(R_H)\rightarrow\mathrm{HH}^*(R_{[^gH]});$$
which can be viewed as: "restriction", "transfer" and "conjugation" maps.

In the main theorem of this short note we verify that the quadruple $$(\mathrm{HH}^*(kH),r^G_H,t^G_H,c_{g,H})_{H\leq G,g\in G}$$ is a Mackey functor, see \cite[\S 53]{Coconet-Todea-The}.
\begin{thm}\label{Coconet-Todea-thmmain} Let $K\leq H\leq G$ and $g,h\in G$. The following statements hold.
\begin{itemize}
\item[i)] $r^H_K\circ r^G_H=r^G_K,~t^G_H\circ t^H_K=t^G_K;$
\item[ii)] $r^H_H=t^H_H=id_{\mathrm{HH}^*(R_H)};$
\item[iii)] $c_{gh,H}=c_{g,^hH}\circ c_{h,H};$
\item[iv)] $c_{h,H}=id_{\mathrm{HH}^*(R_H)}$ if $h\in H$;
\item[v)] $c_{g,K}\circ r^H_K=r^{^gH}_{^gK}\circ c_{g,H}$ and $c_{g,H}\circ t^H_K=t^{{}^gH}_{{}^gK}\circ c_{g,K};$
\item[vi)] $r^G_K\circ t^G_H=\sum_{g\in[K\backslash G/ H]}t_{K\cap {^g}H}^K\circ r^{{}^gH}_{K\cap {}^gH}\circ c_{g,H,~~}$ where $[K\backslash G/H]$ is a system of representatives of double cosets $KgH$ with $g\in G$.
\end{itemize}
\end{thm}

\section{The proof of Theorem \ref{Coconet-Todea-thmmain}}
In order to prove the main result, we need the following lemma. Although the results from the next lemma are easily checkable, for completeness we give the entire   proof. We are inspired by the methods used  in \cite[Lemma 2.1]{Coconet-Todea-Boi}.
\begin{lem}\label{Coconet-Todea-lema} Let $K,H\leq G$ and $g,h\in G$. The following statements hold.
\begin{itemize}
\item[a)] $M$, $N$ and $P$ are projective as left and as right modules.
\item[b)] $R_K\otimes_{R_{[K\cap{}^gH]}}R_{[gH]}\cong R_{[KgH]}$ as $R_K-R_H$-bimodules.
\item[c)] $R_{[g(^hH)]}\otimes_{R_{[^hH]}}R_{[hH]}\cong R_{[ghH]}$ as $R_{[^{gh}H]}-R_H$-bimodules.
\end{itemize}
\end{lem}
\emph{Proof.}

 For the proof of $a)$ we easily see that the first two bimodules are projective  left, respectively right modules since $R_H$ is a parabolic subalgebra of $R_G$, see \cite[Definition 5.1, (Pa2), (Pa1')]{Coconet-Todea-BrHig}.  Since $R_{[gH]}$ is a direct summand of $R_G$ as a right $R_H$-module, and $R_G$ is a projective right $R_H$-module (we use again the argument that $R_H$ is parabolic) we obtain that $R_{[gH]}$ is a projective  right $R_H$-module. Similar arguments show that $R_{[gH]}$ is a projective  left $R_{[^gH]}$-module.

We prove assertion $c).$ The map $$\Phi:R_{[g(^hH)]}\otimes_{R_{[{}^hH]}}R_{[hH]}\rightarrow R_{[ghH]}$$ defined  by  $\Phi(r_1\otimes r_2)=r_1r_2$ for any $r_1\otimes r_2\in R_{[g(^hH)]}\otimes_{R_{[{}^hH]}}R_{[hH]}$ is a well-defined
$R_{[^{gh}H]}-R_H$-bimodule homomorphism. The inverse of this map  is
$$\Psi:R_{[ghH]}\rightarrow R_{[g(^hH)]}\otimes_{R_{[^hH]}}R_{[hH]},$$ given by  $\Psi(r)=\sum_{i\in I} a_i\otimes b_ir$ for any $r\in R_{[ghH]}$, where $1=\sum_{i\in I} a_ib_i,$ and all $a_i\in R_g,$ $b_i\in R_{g^{-1}},$ since $R_1=R_g\cdot R_{g^{-1}}$; here $I$ is a finite set of indices.
Consider another decomposition $1=\sum_{j\in J} a_j b_j,$ where $a_j\in R_g,$ $b_j\in R_{g^{-1}}$ for any $j\in J$ and $J$ is a finite set of indices. Then the equalities
\begin{align*}
\sum_{i\in I} a_i\otimes b_ir&=\sum_{i\in I} 1a_i\otimes b_ir=\sum_{i\in I} \sum_{j\in J} a_jb_ja_i\otimes b_ir\\
&=\sum_{j\in J} a_j\otimes b_j\sum_{i\in I} a_ib_i r =\sum_{j\in J} a_j\otimes b_jr
\end{align*}
show that the last mentioned map is well-defined.

At last, for the assertion b) we fix $t\in K$ and $z\in H$ and we observe that
$$R_t\otimes _{R_1}R_{gz}\cong R_{tgz}$$ as $R_1-R_1$-bimodules. Indeed, this can similarly be shown as statement $c),$ the map given by multiplication has an inverse  that is defined using the relation $1=\sum_{i\in I} a_ib_i,$ where for all $i\in I$ we have $a_i\in R_t$ and $b_i\in R_{t^{-1}}.$
One can observe that, by linearity, the multiplication map extends to a $R_K-R_H$-bimodule isomorphism between
$R_K\otimes_{R_{[K\cap{}^gH]}}R_{[gH]}$ and $R_{[KgH]}.$

\emph{Proof.}\emph{(of Theorem \ref{Coconet-Todea-thmmain}):}

ii) and iv) are an easy exercise if we use \cite[Proposition 2.7 (4)]{Coconet-Todea-LiZh}. Statement v) is  similar to iii) and is left for the reader. Also, the second part of (i) is analogous to the first part and is left as an exercise. The rest of the proof is a consequence of \cite[Proposition 2.7 (1)]{Coconet-Todea-LiZh} and of the bimodule isomorphisms from Lemma \ref{Coconet-Todea-lema}.
Let $X:=R_H,$ viewed as an $R_K-R_H$-bimodule, $Y:=R_{[ghH]}$ viewed as an $R_{[^{gh}H]}-R_H$-bimodule, $Z:=R_{[g(^hH)]}$ viewed as an $R_{[^{gh}H]}-R_{[^hH]}$-bimodule and $U:=R_{[hH]}$ viewed as an $R_{[^hH]}-R_H$-bimodule.
For (i) we have
$$r^H_K\circ r^G_H=t_{X}\circ t_{M}=t_{X\otimes_{R_H}M}=r^G_K.$$
For iii) we have
\begin{align*}
&c_{gh,H}=t_{Y}
\mbox{ and }\\
&c_{g,^hH}\circ c_{h,H}=t_{Z}\circ t_{U}=t_{Z\otimes_{R_{[^hH]}}U};
\end{align*}
From Lemma \ref{Coconet-Todea-lema}, c) we know that  $Y\cong Z\otimes_{R_{[^hH]}}U$ as $R_{[^{gh}H]}-R_H$-bimodules, hence iii). Consider the $R_K-R_G$-bimodule $A:=R_G$, the  $R_K-R_H$-bimodule $B:=R_G$, the $R_K-R_{[K\cap^gH]}$-bimodule $C:=R_K$ and the  $R_K-R_H$-bimodule $D:=R_{[KgH]}$.  For vi) we have
$$r^G_K\circ t^G_H=t_{A}\circ t_{M}=t_{A\otimes_{R_G}M}=t_B,$$
since $A\otimes_{R_G}M\cong B$.

By \cite[Proposition 2.11 (i), 2.12 (iii)]{Coconet-Todea-LiTr} and since $C\otimes_{R_{[K\cap^gH]}}P\cong D$ as $R_K-R_H$-bimodules (see Lemma \ref{Coconet-Todea-lema}, b))   the following equalities hold
\begin{align*}
\sum_{g\in[K\backslash G/H]}t_{K\cap^gH}^K\circ r_{K\cap^gH}^{^gH}\circ c_{g,H}
&=\sum_{g\in[K\backslash G/H]}t_{C\otimes_{R_{[K\cap^gH]}}R_{[{}^gH]}\otimes_{R[{}^gH]}P}\\&=\sum_{g\in[K\backslash G/H]}t_{C\otimes_{R_{[K\cap^gH]}}P}\\&=\sum_{g\in[K\backslash G/H]}t_{D}\\&=t_{B}.
\end{align*}

%%%BIBLIOGRAPHY

\end{document}